**Ola Bratteli and his diagrams**

**Tone Bratteli, Trond Digernes, George Elliott, David E. Evans, Palle E. Jorgensen, Aki Kishimoto, Magnus B. Landstad, Derek Robinson, Erling Størmer**

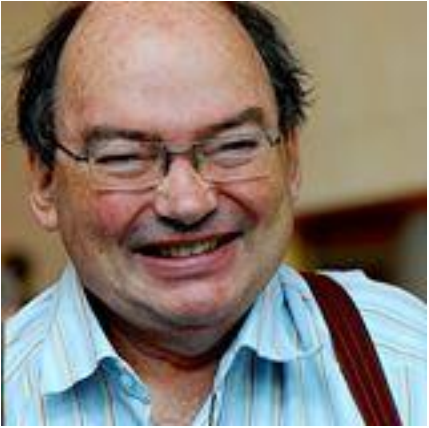

### Short CV

Born October 24, 1946, died February 8, 2015.

Graduated with distinction from the University of Oslo in 1971 and took his doctorate there in May 1974. He was a research fellow at New York University 1971-73, had various post.doc-positions 1973-77, he was an associate professor at the University of Oslo 1978 -79, a full professor at the University of Trondheim (now NTNU) 1980–91, and from 1991 at the University of Oslo.

Ola's father Trygve Bratteli was a Norwegian politician from the Labour Party and Prime Minister of Norway in 1971–1972 and 1973–1976. During the Nazi invasion of Norway, he was arrested in 1942, was a Nacht und Nebel prisoner in various German concentration camps from 1943 to 1945, but miraculously survived. Ola's mother Randi Bratteli was a respected journalist and author of several books.

Ola is survived by his wife Rungnapa (Wasana) and their son Kitidet.

According to MathSciNet Ola Bratteli has 113 publications with 21 coauthors, he received various awards, and was a member of the AMS for 43 years.

For a detailed biography, see MacTutor History of Mathematics Archive:
http://www-history.mcs.st-and.ac.uk/Biographies/Bratteli.html

### Bratteli diagrams

### George Elliott

One of Bratteli's most important discoveries, I think, was what is now called a Bratteli diagram, which he found as a way of codifying the data of an approximately finite-dimensional (AF) C*-algebra (the completion of an increasing sequence of finite-dimensional C*-algebras -- finite direct sums of matrix algebras), in a far-reaching extension of the thesis of Glimm (who considered sequences of simple finite-dimensional C*-algebras, single full matrix algebras, with unital embeddings -- the non-unital case was later studied by Dixmier).

In one sense, Bratteli diagrams had perhaps already been invented, as, for one thing, the idea is so simple -- a sequence of horizontal rows of points, with numbered lines connecting points of each row to the points of the row below, these numbers recording the multiplicities of the partial embeddings of the simple direct summands at one stage of the sequence of algebras into those at the next stage. (If, then, numbers are assigned to the points in the first row, which may be assumed to be a single point and the assigned number just one, numbers are then automatically assigned by the multiplicities to the points in later rows -- this assumes that one is talking about unital C*-algebras, but this involves no loss of generality.)

Pascal's triangle is a Bratteli diagram -- with the multiplicities equal to one, and the numbers appearing in the rows being of course the successive degrees of binomial coefficients. The AF C*-algebra this diagram encodes arises in physics as the gauge-invariant subalgebra of the C*-algebra (also AF) of the canonical anti-commutation relations!

Bratteli was not content just to look at the diagrams -- he isolated the equivalence relation between them which is determined by isomorphism of the associated C*-algebras. This was prophetic, as he in fact was noticing that the diagrams formed a category, in which his equivalence is just isomorphism. It was later noticed that this category is equivalent to the category of ordered groups arising from the algebras in question via the $K_0$-functor. This led eventually to a K-theoretical classification of an enormous -- still evolving -- class of (what are called) amenable C*-algebras -- analogous to the classification by Connes and Haagerup (and others) of amenable von Neumann algebras.

Bratteli diagrams arose in a fundamental way in Jones's theory of subfactors. Given a subfactor of Jones index less than four, the increasing sequence of relative commutants in the Jones tower are finite-dimensional, and so give rise to an AF algebra. Its Bratteli diagram is periodic with period two, with the step from the second row to the third being just the reflection of the step from the first row to the second. Both steps are obtained from a Coxeter-Dynkin diagram, by dividing its vertices alternately into a first row and a second row, so that the edges (already equipped with multiplicities!) connect the two rows and constitute the first step of a Bratteli diagram---the second step of which is obtained by interchanging the two rows.

Bratteli diagrams with an order structure were introduced by Vershik to describe a measurable transformation, and were adapted by Herman, Putnam, and Skau to describe a minimal transformation of the Cantor set. Using this description, Giordano, Putnam, and Skau classified the orbit structures of such transformations, the invariant being (in the generic case) the ordered $K_0$-group of the associated C*-algebra, also classified by this ordered group -- no longer an AF algebra, but belonging to a class determined by $K_0$ plus $K_1$, the latter group being in the present case just the integers.

Below in Figure 1 is a Bratteli diagram representing the so-called GICAR-algebra. GICAR stands for gauge invariant CAR (which again stands for canonical anticommutation relations). Note the resemblance to Pascal's Triangle.

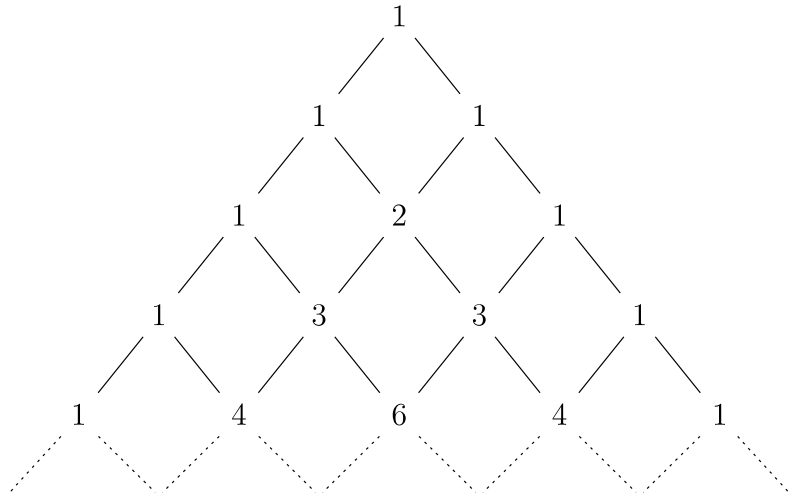

The corresponding inductive chain system (depicted vertically) is:

$$\mathbb{C}$$
$$\downarrow \phi_0$$
$$\mathbb{C} \oplus \mathbb{C}$$
$$\downarrow \phi_1$$
$$\mathbb{C} \oplus M_2(\mathbb{C}) \oplus \mathbb{C}$$
$$\downarrow \phi_2$$
$$\mathbb{C} \oplus M_3(\mathbb{C}) \oplus M_3(\mathbb{C}) \oplus \mathbb{C}$$
$$\downarrow \phi_3$$
$$\mathbb{C} \oplus M_4(\mathbb{C}) \oplus M_6(\mathbb{C}) \oplus M_4(\mathbb{C}) \oplus \mathbb{C}$$
$$\downarrow \phi_4$$
$$\vdots$$

where the (injective) connecting homomorphisms are given by

$$\phi_0(a) = a \oplus a,$$
$$\phi_1(a,b) = a \oplus (\begin{smallmatrix} a \\ & b \end{smallmatrix}) \oplus b,$$
$$\phi_2(a,B,c) = a \oplus (\begin{smallmatrix} a \\ & B \end{smallmatrix}) \oplus (\begin{smallmatrix} B \\ & c \end{smallmatrix}) \oplus c,$$
$$\phi_3(a,B,C,d) = a \oplus (\begin{smallmatrix} a \\ & B \end{smallmatrix}) \oplus (\begin{smallmatrix} B \\ & C \end{smallmatrix}) \oplus (\begin{smallmatrix} C \\ & d \end{smallmatrix}) \oplus d,$$
$$\vdots$$

where blank means one or more 0's, a, b, c are complex numbers and B, C are matrices.

Ola was always a fountain, or mountain, of good sense. Once, when he was visiting Toronto, we went to the gym to go swimming. I was a member, but it was a little murky what guest privileges I had. I thought I had them, as a faculty member, and proceeded to explain that. This met with resistance, which gradually became more protracted. In the meantime, Ola slipped past the desk, picked up a towel, and half an hour later had finished his swim--- perhaps even had a sauna, too. At that point I gave up and we left.

## Ola – a child of peace, a man of the outdoors, and a family man

### Tone Bratteli

On United Nations Day, October 24, 1946, a child of peace was born in a crowded maternity department in Oslo. He was one of many in the baby boom that followed World War II.

This event was no matter of course. My father came home in 1945 after years in extermination camps in Germany. He survived by a hair's breadth. Soon after his return he met my mother. Her father had also come home from concentration camps. The two found each other quickly, despite my father's shyness. My mother's sociable nature made up for that. And in October 1946 he could lift his son up in the air in joy. He had not been sure whether he would be able to have children after the appalling treatment in the camps.

A year and a half later – May 8, 1948 – I was born. And after another three years – May 20, 1951 – our little sister Marianne arrived. Ola's sometimes insistent little sister.

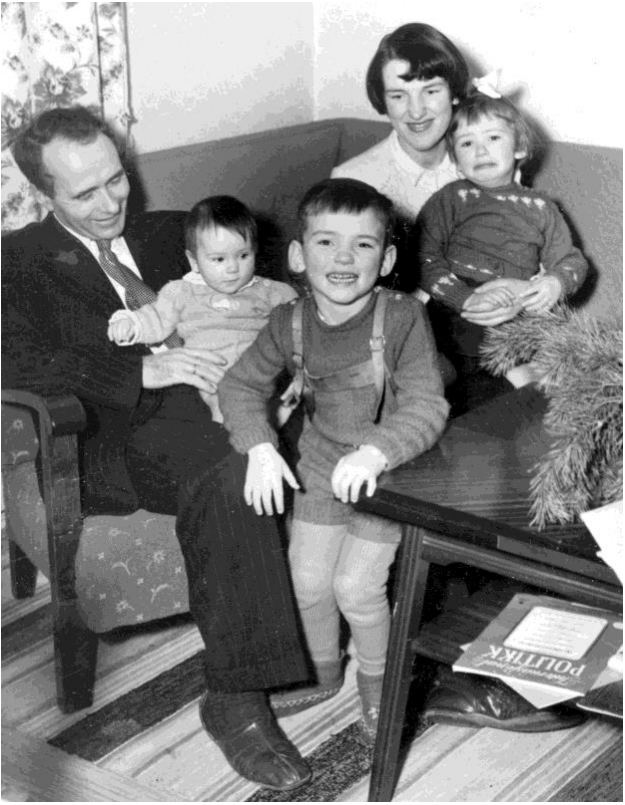

My father was on his way into politics, with the result that – when Ola was five years old – we moved to a so-called official residence in which our parents could also carry out social duties. On the day we moved, Ola and I scrambled around searching for our old home. The little we had with us was swallowed up by the huge rooms. Many leading politicians from other countries came to this apartment. We children hid away in our rooms and were not especially eager to introduce ourselves. Ola was a quiet boy, but one evening he came in to join the house guests and started tugging on the men's ties. This was a big surprise. What had happened to Ola? He had been in the kitchen, where some half-empty wine glasses had caught his eye. After tasting that juice, the shy, self-conscious boy was briefly transformed into a party animal.

But Ola was soon ready for Bolteløkka School and a meeting with the subject that became his passion – mathematics.

My mother was rather taken aback the day she discovered that Ola had wallpapered his room with equations. He did not play football or hang out with friends – he solved equations and went for very long walks in the forests and the mountains. He also did my maths homework for me, because to me that subject was a struggle. Ola achieved the best grade in mathematics, mine was the worst.

As a young man, Ola spent many hours skiing, and the trips could easily reach 50 to 60 kilometres. He also tried to find detours to make them even longer. Loyal skiing friends accompanied him on some of the trips. When he came home, a slice of bread or two was not enough. He ate the whole loaf.

We had a summer cabin at Nesodden, half an hour from Oslo by boat. We spent our childhood summers there. At Nesodden we learned to swim, gather blueberries and mushrooms, empty the outdoor toilet, play badminton and find spooky ghost trails in the surroundings. Those were good summers. A time of freedom and peace of mind. When Ola was a little older, he sometimes cycled to the cabin instead of taking the boat. It was a 50-kilometre journey.

Early on, Ola showed an interest in music and visual arts as well. He made an attempt to teach himself to play the cello. It was not a success. But he spent a lot of time listening to music and he went to exhibitions.

Our sister Marianne is an artist, so my brother and sister had something in common there. Later in life, when I travelled to several continents in my job as a journalist and in other connections, Ola always knew which gallery I should head for if there was only time for one. The composition of a picture and the solution to a mathematical puzzle must have something in common.

Ola graduated from the University of Oslo in 1971. He was awarded the best possible grade, and was what we call "reported to the King". That meant that the entire government was informed about his academic accomplishments. It was a big day for everyone in the family, but perhaps most of all for my father. He was the Prime Minister who came from such a poor background that he never completed his education himself, but was self-taught. Now he was to inform the King and government about his son's academic triumph. As usual, Ola was unassuming and self-conscious, but he was no doubt satisfied.

Now a whole world of academic challenges lay before him. It was difficult to be completely anonymous and work undisturbed with his own research in Oslo. Being the child of a Prime Minister has many sides. You have to take a variety of comments and media coverage in your stride. When Ola went to New York in 1971, it was a kind of escape.

From his apartment in Greenwich Village, it was easy to get to theatres and exhibitions, and Ola soaked up everything he could get to. One day when it snowed in New York – that does happen, after all – he skied up and down Fifth Avenue. Finally, skiing weather in the Big Apple...

Ola returned home in 1973 and took his doctorate in 1974. This was during Dad's second term as Prime Minister, and also led to coverage in the media.

In his personal life, Ola was simply a very kind, generous, caring and family-loving man. With a twinkle in his eye, it was easy for him to establish rapport with children.

Ola was interested in genealogy. He collected a great deal of material about the family's deep roots in Norway. They could be traced back to the 13th century. He was especially interested in his paternal heritage.

Ola cared about all his aunts, uncles, and cousins. Our father was one of 11 children and our mother one of 3, so there were quite a few of us. Sometimes Ola would cycle to a cousin living in Kongsvinger to discuss our family tree. The trip to Oslo and back is 186 kilometres, and it was no problem for him to complete it in a day. After all, he was on a quest for our identity.

During the years he lived abroad, he kept in touch with his mother by phone and wanted to know how we were doing. While he was in Australia, he did not always come home for Christmas. A barbecue on the beach in the sunshine was more tempting than snow, sleet and freezing temperatures. Perhaps he had covered enough distance on skis.

Before Ola moved back to Norway, he was fortunate enough to meet Rungnapa (Wasana). He was to share almost half his life with her.  Kitidet – his son – was his pride and joy.

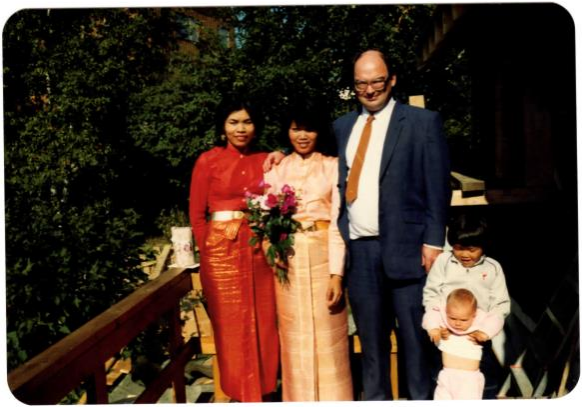

Trondheim 1986

Rungnapa and Ola travelled widely all over the world. There were holidays, but she also accompanied him to conferences and on visits to universities. After a while, they settled down in Norway. They enjoyed good years together.

In the last years of his life, Ola's health deteriorated. Rungnapa was enormously supportive and did all she could to make it possible for him to live at home as long as possible. He was a lucky man.

Now Rungnapa and Kitidet have moved back to Phitsanulok in Thailand. Ola and Rungnapa built a house in her home city several years ago and had no doubt planned to spend the winters there eventually. It did not turn out that way. Ola died so early. But outside the house there is a small temple for Ola. So in a way he is there too.

## Trond Digernes

I first met Ola around 1970 when we were both master students at the University of Oslo. At that time Ola was a slender young man with a passion for the outdoors, especially mountain hiking. He also had an adventurous side: he had spent a summer on the arctic archipelago of Svalbard as a geologist's assistant.

During the year 1970-71 there were three of us who spent a lot of time together: Ola, John Erik Fornæss and myself. We played bridge, took a skiing vacation in the Norwegian mountains, and went mountain hiking in the summer. At the end of summer 1971 our roads parted. We all went to the US for PhD studies, but to different institutions: Ola to Courant Institute, John Erik to the University of Washington, and myself to UCLA.

The next summer Ola came over to the West Coast along with two other Norwegians. The four of us went hiking in the Sierra Nevadas for several days, carrying tents, sleeping bags, and provisions on our backs. - My next adventure with Ola was in the summer of 1973 when we met in Iceland, along with my brother. We rented a four-wheel drive, drove across the roadless interior of the island, got stuck in rivers, and slept out in the open. After Iceland Ola returned to Courant whereas I was headed for a year's stay at CNRS, Marseille.

At CPT/CNRS, Marseille, the year 1973-74 was organized as a special year dedicated to operator algebras and mathematical physics. It attracted several high-powered researchers – among them Masamichi Takesaki and Alain Connes – and people like Derek Robinson and Daniel Kastler were already there. To Ola this sounded like a good place to be, so he left Courant and joined the Marseille group in January 1974. This was also when he started his long-time collaboration with Derek Robinson.

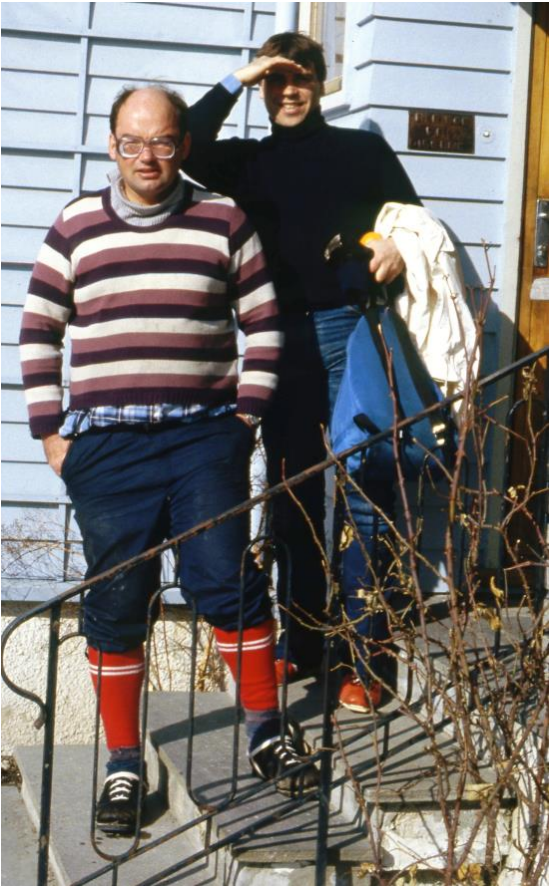
Ola and Trond, 1983

After Marseille I returned to UCLA to finish my PhD, and subsequently went on to Trondheim for a position at the University of Trondheim (later reorganized as The Norwegian University of Science and Technology). Ola held various positions in the period 1974-80, until he came to the University of Trondheim as a professor in the fall of 1980. During the decade 1980-90 I joined Ola and Derek on several occasions for discussions at ANU, Canberra, and Derek also visited Trondheim. This resulted in a few joint publications, sometimes also with other co-authors. In 1991 Ola became professor at the University of Oslo, a position he held until his retirement.

I was only involved in a fraction of Ola's mathematical work, but since we spent much time together, we had many interesting discussions over the years. Given Ola's incisive mind and deep understanding of everything he was involved in, it was always a rewarding experience to discuss with him. He is dearly missed, both as friend and colleague.

## Ola Bratteli, friend and mathematician

### Erling Størmer

I met Ola for the first time when he was about to start on his Masters thesis. Then he was a dark-haired lad with a beard, radiating health and fitness, who often went for very long skiing trips. It soon struck me that he was a highly effective person who could pick up new theory extremely quickly. By then he had taught himself a great deal about the field of operator algebras, which he wanted to work on. This branch of mathematics started in the years around 1930 when people wanted to develop mathematical theory for quantum mechanics. The theory was developed further by mathematicians including John von Neumann, but there were few active participants until the late 1960s. Then the physicists joined in, and operator algebras became a popular field. So Ola's timing was excellent when he passed Masters examination in 1971 with top grades. When it became clear to me how good the thesis was, I said to Ola that we had made a big mistake; this thesis should have been used for a doctoral degree. And it was precisely the results here – see [1] – that made Ola well-known as a mathematician from an early stage.

In 1959, James Glimm published a PhD thesis on operator algebras, which were achieved by taking an infinite union of an increasing family of n x n matrices. This became a famous piece of work, and when Ola was to write his thesis, I suggested that he could extend Glimm's results to some more general operator algebras. Ola quickly discovered that he could generalize Glimm's work by studying infinitely increasing unions of matrix algebras. He then ended up with an infinitely large diagram that described all the inclusions, which thus also described how the operator algebra was constructed. His main finding was that this diagram fully described the operator algebra, enabling a classification of all such operator algebras, now known as AF algebras, an abbreviation for "approximately finite dimensional operator algebras". This result proved far more important than Ola and I had anticipated. The reason was that analogous diagrams could be used in other areas of mathematics, with very good results. They are now called Bratteli diagrams, and are a well-known concept among mathematicians.

After graduating, Ola studied for two years at New York University, where Glimm was, and there was an excellent environment in mathematical physics. As well as his studies, Ola spent a great deal of time in New York on cultural life and on eating well. During those two years, his appearance changed in some ways; much of his hair disappeared, his beard was gone, and he put on so much weight that when I met him a couple of years later, I did not recognize him, so I introduced myself to him.

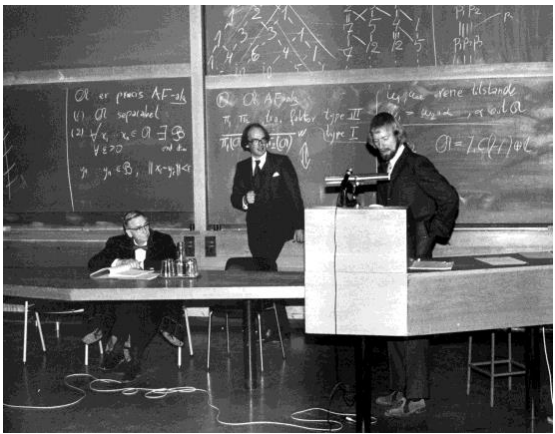

Ola's Ph.D defense, from left: The Dean, Ola, Gert K. Pedersen.

Ola returned home and took his doctorate in Oslo in 1974. After that, he went abroad again. At that time, there was a very active and high-quality environment in Marseilles led by physicists who developed the theory of quantum physics in operator algebras. Ola spent a great deal of time and became a very good friend and research partner with an English colleague, Derek Robinson. The two soon published several works together, and they continued with this for the rest of Ola's career. Robinson was the co-author of Ola's last two research articles in 2008. They also wrote a two-volume work together on operator algebras; see [3] and [4]. Intended especially for physicists, it became very popular, and has become a standard reference for physicists in the field of operator algebras.

After a few years, Robinson moved to Australia, so Ola travelled there several times. On one of the journeys he stopped in Thailand, where he was fortunate enough to meet Wasana, whom he married. Together, they had a good life. But Robinson was not Ola's only close research partner. He also did a great deal of work with four other colleagues, who became his close friends. They were George Elliott, David Evans, Palle Jørgensen and Akitaki Kishimoto.

In the years after 1980, Ola published 84 works, of which he was the sole author of only 4, but with fully 74 in which at least one of those I have mentioned was a co-author. Often, there were other co-authors as well. Ola was altogether a very popular person to work together with. Among other things, as I have mentioned, he was very effective when he worked. Everyone liked Ola. When there were several of us together, he was not a man of many words, but in private he would talk. Ola radiated a good spirit, and it was easy to become fond of him. He was a person with a warm heart.

As our colleague Sergio Doplicher in Rome expressed it so well when he heard about Ola's death: "We will greatly miss Ola, as a scientist and as a friend, his kind and soft approach to the others, always based on understatement and totally antipodal to any form of arrogance, but with the constant dressing of his subtle humour, will stay with us for the rest of our lives."

Ola was a professor at the University of Trondheim during the period 1980–1991, and he was Head of Department during the period 1981–1984. From 1991, he was a professor at the University of Oslo. As a lecturer, he was perhaps best suited to teaching at master's level. His style fitted in well there, and he gave excellent lectures that were very clear and well prepared.

Ola Bratteli passed away at the age of 68 after several years of steadily declining health. It was an extremely sad experience to see how he had become weaker each time I saw him in recent years. His strength began to fail at a fairly early stage. His last research articles were published in 2008, and after that he had little energy to do more. So his brilliant career as a mathematician came to an end far too early. Ola will be remembered and missed for a long time, both as a mathematician and as the wonderful person he was.

## Life with Ola

### Derek Robinson

I met Ola at the beginning of 1974 in Marseille. He was introduced to me by Trond Digernes. It was to be a pivotal moment in each of our lives, although we had no premonition of this at the time. There were many unpredictable consequences: a few years later Trond was to meet his wife Hallie at an open air opera in Sydney, Ola became the owner of a mushroom farm in Northern Thailand and Ola and I were to write a book that is still bought, read and regularly cited 35 years later.

At that time I was Professor of Physics at the Luminy Campus of the University of Marseille. Ola and Trond were graduate students in mathematics; Ola in Oslo with Erling Størmer and Trond in Los Angeles with Masimichi Takesaki. Luminy, largely under the influence of the late Daniel Kastler, then had developed a strong visitor's program in mathematical physics and this explained Ola and Trond's presence.

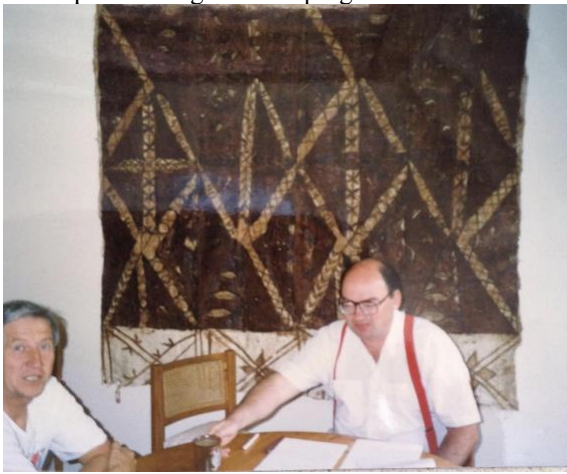

## Derek and Ola

My collaboration with Ola began the day we met. I explained to him some of the ideas I had about quantum dynamics and derivations on C*-algebras and shortly after we wrote our first paper on these topics. By 1976 we had coauthored three other papers. Although, in principle, we were both in Marseille for this period in practice we were both often travelling and collaborated by mail. This was to set a pattern for the next 40 years.

Although Ola was still a graduate student he had already completed his now famous, influential and often cited work on AF-algebras. I, on the other hand, was involved with several different areas of mathematical physics. When we met we had different interests, different backgrounds and quite different personalities. Ola was quiet, well organized and introspective, characteristics I did not share. Ola rapidly assimilated the Mediterranean lifestyle, the sun, the sand and the seafood, especially the seafood. It was a time of calm cooperation amid the chaos of French academic life.

Our collaboration took an unexpected turn in autumn 1975 when I received an invitation to lecture at a Mathematics Research Institute in Australia. It was during this trip that I began to think about writing a book on operator algebras and their applications in physics. I thought that Ola and I, with our disparate backgrounds and successful working relation, would be able to do justice to the subject and its recent developments. On my return to Marseille I discussed the idea with Ola and was pleased that he did not dismiss it immediately.

Initially the book was intended as a relatively short term project; a monograph of 3–400 pages with the early chapters on mathematical background and the later chapters on applications to quantum statistical mechanics. We began in September 76 with a chapter on the basic theory of operator algebras and quickly realized that we would exceed the estimated length. By September 77 we had the equivalent of 500 printed pages of material and had not reached the chapters on applications. So the short term project turned out to be a long term project and the book changed from one volume to two. Then a second effect of my visit to Australia complicated things. I decided to move from Marseille to Sydney. This meant that the second volume was largely written with Ola in the Northern Hemisphere and me in the Southern Hemisphere. Fortunately the writing procedure we adopted was quite robust.

We began each chapter with a tentative sketch of the intended sections. Then, after discussing the general presentation of the material in each section, we began to draft alternate sections. For example, I drafted the first section on general algebraic structure and Ola drafted the second section on representations of algebras. Then we exchanged drafts and edited each other's work. This process would be repeated until we were each satisfied with the outcome. I am not sure whether this is a standard procedure but it worked well with us. The editing was not a superficial process. We often had different notions of the relative significance of the material and the emphasis to be given to various statements and results. At times my first draft would be completely changed by Ola and vice versa. Somehow the process always reached equilibrium after a reasonably short time, with one exception the section, on modular theory. This took seven exchanges before we were both satisfied. This procedure had various advantages. It naturally introduced a uniformity of style. It also gave a fairly foolproof method of avoiding error, although we were not completely successful in that respect.

The first volume of the book was completed by September 77 which left three months to complete the second volume of the book before my departure with the family for Australia. In that time we managed to write about 40% of Volume 2, aided in part by discussions with Akitaka Kishimoto who had just arrived in Marseille as a postdoc. My recollections of the first half of 1978 are rather blurred. I was fully occupied settling into a new environment and adapting to the Australian academic framework. I think the book writing came to a halt although we did proofread Volume 1. The next time I saw Ola was in the Australian midyear break, June 78. I returned to Marseille for a month and our collaboration started up again. Akitaka was still in Marseille and the three of us collaborated on material which was incorporated into Volume 2. About the time I returned to Australia Ola returned to Oslo and Akitaka moved to Ottawa.

The next eleven months Ola and I continued our collaboration regularly by mail. There was a reliable airmail system between Australia and Norway. Typically one of us would write a section of manuscript mail it to the other who would edit it and mail it back. The whole process took two weeks. It was not ideal but we adapted to the routine. Then, however, there was a disaster. The Sydney international mail exchange was closed by a strike which lasted two months. Again our collaboration ground to a halt.

After the Sydney postal strike ended we continued on the final chapter of the book. Then in June 79 we both returned to Marseille and began the final work. We shared an apartment in Bandol on the coast East of Marseille and worked non-stop. We had four weeks to complete the book and it took us three and a half. The second volume finally appeared in 1981. So the total operation of writing and publishing the 1000 page two volume book took about three years. That was not the end, however, we returned to it again twice preparing the second edition, but that is a different story.

After the book was finished collaboration with Ola and Akitaka continued. They both visited Sydney, and subsequently Canberra, for extended periods. Ola was lucky to survive his first visit. Whilst touring the almost deserted roads of a National Park he reverted, European style, to drive on the right side, the wrong side for Australia. This led to a head on collision which wrecked both cars, fortunately without any personal injuries.

Ola realized that by timing his Australian visits correctly he could ensure that it was always summer. It also had the advantage that we could spend a maximum of time working at the family beach house. Ola would then take his daily swim before heading to the local oyster shop. He bought freshly harvested oysters by the bag to be opened, seasoned with a lemon from the garden, and savored as the sun went down.

'Those were the days my friend, I thought they'd never end'

## A tribute to Ola Bratteli

## Aki Kishimoto

Before I first met him in September of 1977, I must have read his early paper [1] of 1972 and his more recent series of papers on unbounded derivations, mostly written with D.W. Robinson, which would be incorporated into their forthcoming books [3,4], because his image had been firmly established in my mind as a formidable mathematician I could hardly be compared with. And he was; and I think I was quite lucky to come to know him in the earliest possible days, which brought me chances to collaborate with him ever.

Ola and Derek were writing the book in Marseille, a kind of book I love and could get familiar with, when I visited in 77. Derek soon left for Sydney (where we would meet later) but Ola and I spent almost a year together there. Though he was just one year my senior he was a kind of mentor in mathematics and everything else during the stay. We often ate out dinner and sometimes spent weekends together. He also arranged me to attend two conferences, one in Naples and the other in Oslo in 1978, sort of encouraging me against my timidity. In all these activities (except for the Oslo) when transportation was required Ola was in charge, which eliminated a practical worry from me.

In one outing I remember we hiked on rocky coastal paths, where I learned to mumble 'Bonjour' to strangers we encountered. (Later I found this practice ubiquitous and Ola would say 'Konnichiwa' awkwardly when we hiked in Japan.) We were not a type of enjoying sports, but he was better than me in everything. So when we had to descend a steep slope made of gravels at one point, I was surprised to find myself rather enjoying skidding down while Ola tried to walk down steadily, as if reflecting his meticulous style of doing mathematics.

Ola's paper [1] on AF algebras became an inspirational source for me. This class might look rather special but now we might say *if a C\*-algebra is not obviously not AF then it would be AF*. As a touchstone of this credo Ola and I examined the fixed point algebra $F_\theta$ of a $C^*$-algebra $C_\theta(u,v)$ generated by two unitaries $u,v$ with $uv = e^{2\pi i\theta}vu$ with $\theta$ irrational, under the period two automorphism σ; $\sigma(u) = u^{-1}, \sigma(v) = v^{-1}$. ($C_0(u,v)$ with $\theta$=0 still can be defined and is the continuous functions on the torus while $F_0$ is on the sphere or rather a pillow with four corner points.) We managed to show $F_\theta$, a non-commutative pillow, is AF (when $\theta$ is irrational). When Ola was finishing the paper assiduously as always, I had to leave him alone in Sapporo. A note I found the MS with on return gave as an excuse of a messy proof a missing hole he stumbled onto. ???? (This was published as Non-commutative spheres III in J. Funct. Anal. 147 (1992).) Another notable result in [1] was that any two irreducible representations of a simple AF algebra are bridged by an automorphism, which we established in a different setting, published as *Homogeneity of the pure state space of the Cuntz algebra* in J. Funct.Anal. 171 (2000). This gave me a hope to attack a more general case and I did with some help from others.

Ola had many works; not all I can cope with. Among them the books [3,4] is a good reference for those in the field of mathematical physics, me included. Our last work *Approximately inner derivations*, Math. Scand. 103 (2008) with Derek is an attempt to shed light on a topic dealt there, which yet haunts me to this day.

### Palle E . Jorgensen

Ola Bratteli (1946 – 2015) was a Norwegian mathematician who over a long period, starting with 1972 has had a profound influence on modern analysis; especially themes connected with operator algebras, classification, non-commutative harmonic analysis, and representation theory. While Ola's first paper was solo, almost all that followed were joint. And there were many collaborators. A glance at MathSciNet demonstrates that the papers are widely cited.

In January of 2017, a week-long conference was organized, with the aim of presenting some of the many collaborative advances in mathematics and its applications involving Ola. The conference took place at a Danish-Norwegian conference center, Lysebu, a scenic location up in the mountains overlooking Oslo.

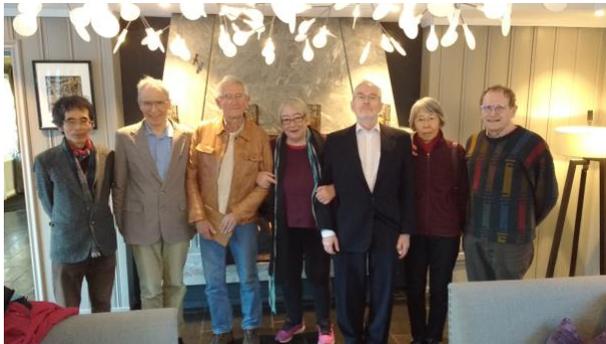

Aki Kishimoto, Palle Jorgensen, Derek Robinson, Tone Bratteli, George Elliott, Reiko Kishimoto, David Evans.

Below is a summary presented by one of Ola's collaborators, Palle Jorgensen. Some of the Bratteli-Jorgensen papers have also included other co-authors, especially Derek Robinson, David Evans, George Elliott, A. Kishimoto, Geoff Price, Fred Goodman, Ki Hang Kim, and Fred Roush; but this is only a partial list. My own collaboration with Ola started by chance, and dates back to the mid 1970ties, and it lasted for four decades. Our early work was in non-commutative geometry, and our later research moved in a diverse number of directions.

While Ola's first and solo paper introduced what is now known as Bratteli-diagrams, the subsequent research moved in a number of different directions; some sketched below. Bratteli diagrams are combinatorial graphs of vertices organized into level, and with edges between vertices at neighboring levels. Their subsequent applications include K-theory, classification of C*-algebras, infinite tensor products, and of representations, and topological dynamics. They are even used as tools in the programming of (finite) fast Fourier transforms. There is a Wiki page, and 20,000 hits with a Google search.

Ola liked to spend extensive periods on research visits in different corners of the World, and his arrangement at University of Oslo fortunately allowed this. Even when Ola was in Oslo, he had grant support for research visitors. My own (Jorgensen) collaborations involved my visiting Ola in Oslo. In addition, we both made research visits to Derek Robinson in Australia, to Dai Evans in the UK, and to George Elliott in Toronto. And there were research visits to workshops at the Mittag-Leffler Institute in Stockholm, Oberwolfach in Germany, and to Banff in Canada. In all, I have 30 joint research publications with Ola, and a book. Two of the papers are in fact AMS Memories. Other papers appeared in such journals as J. Functional Analysis, in *Ergodic Theory Dynam. Systems, in Appl. Comput. Harmon. Anal.*, in *J. Geom. Anal.*, and in *Nonlinearity*.

Early joint research includes the themes: Lie algebras of operators, smooth Lie group actions on noncommutative tori, a study of decomposition of unbounded derivations into invariant and approximately inner parts. These topics are part of a systematic analysis of unbounded ∗-derivations as infinitesimal generators in operator algebras; with direct connections to quantum statistical mechanics. Indeed the motivation for the study of unbounded ∗-derivations is from quantum physics. Hence, the unbounded ∗-derivations serve as generators for the dynamics in particular quantum systems, but realized in the Heisenberg picture; i.e., when the dynamics is a one-parameter group of automorphisms, acting in the particular operator algebra representing the quantum observables. Other applications include non-commutative geometry, as envisioned by Alain Connes.

Our later joint research moved more in the direction of representations and certain applications. During my many research visits to Oslo I came to know Ola's family, and many Norwegian friends. One of Ola's sisters, Marianne is an accomplished artist. She took over a house from Dad, located on an idyllic island in Oslo fjord (Ola's father was Trygve Bratteli, a renown Norwegian prime minister). On weekends, Ola and I were invited to visit the house on the island. I was the Danish visitor. Sadly, Ola's health declined towards the end, but I am happy to have had the benefit of many intense research experiences from the early part of our careers.

My most recent, and substantial, joint work with Ola was wavelets. That part includes a book "Wavelets through a looking glass. The world of the spectrum," which presents the subject from a representation theoretic viewpoint. The research papers leading up to the book include the themes: Wavelet filters and infinite-dimensional unitary groups, compactly supported wavelets and representations of the Cuntz relations; and wavelets from isometries and shifts, and multiresolution wavelet-analysis of scales arising from representations of the Cuntz algebras.

A common theme in my joint research with Ola is my insistence on the central role to be played by representations and decomposition theory. Some of our early work dealt with representations of Lie groups, of C*-algebras, and of multi-scale systems. Often the research was collaboration teams. For example, our work on analysis of Lie groups included Derek Robinson (and others), covering asymptotic laws for periodic subelliptic operators on Lie groups (scaling laws for small time, and for large time); a topic in turn inspired by Hunt processes. This theme also entailed elliptic second order operators and associated unitary representations of Lie groups and Gårding inequalities; and it expanded into a systematic study of heat semigroups and their connection to integrability questions for Lie algebras of unbounded operators.

A quite different representation theoretic theme was the theory of numerical AF-invariants (AF is for approximately finite-dimensional), representations and centralizers of certain states on the Cuntz algebras; and a related but different study of combinatorial notions we called iterated function systems and permutation representations of the Cuntz algebra. Why the Cuntz algebras? As C*-algebras they are simple and purely infinite. They are easily presented as algebras on a finite number of generators ($O_n$ if n is the number of generators), subject to relations, the so called Cuntz relations. Their representations are surprisingly versatile, and they play a key role in our understanding of such multi-level systems as wavelet multi-resolution scales; in addition to multi-band filters in signal processing.

It is known that the equivalence classes of irreducible representations of $O_n$ do not admit a Borel labeling; but nonetheless, in one of the Bratteli-Jorgensen papers we show that the wavelet representations are irreducible, and the equivalence classes are labeled by an infinite-dimensional unitary loop group.

Later joint work between Ola, me, K.-H. Kim and F. Roush was inspired by Bratteli-diagrams. It is known that the non-stationary case defies classification (order-isomorphism is undecidable). But we discovered that the stationary case could be decided by explicit classification numbers and associated finite algorithms. In this work, for the stationary dimension groups, we obtained explicit computation of numerical isomorphism invariants. We proved decidability of the isomorphism problem for stationary AF-algebras and the associated ordered simple dimension groups.

## Ola and Orbifolds

## David E Evans

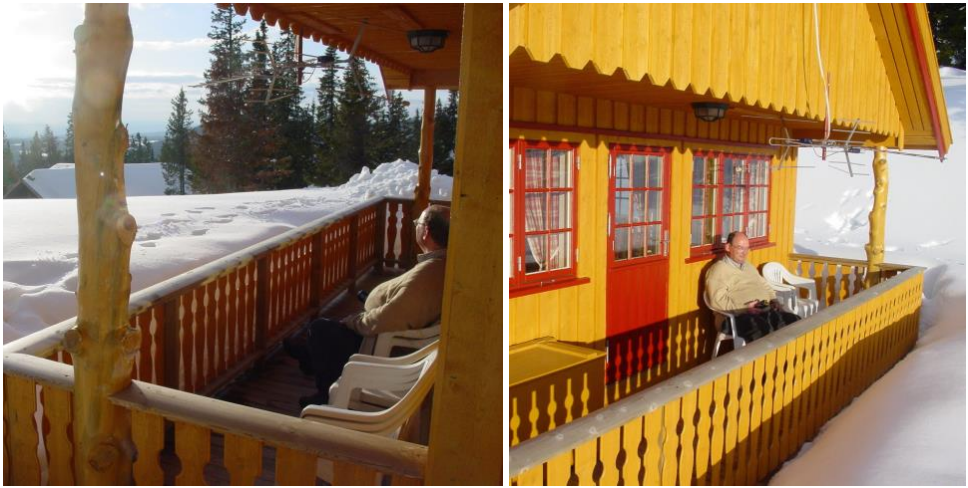

Sjusjøen, Lillehammer, Easter 2004

I first met Ola in the new year of 1977 when I was a postdoc in Oslo. Ola was then based in Marseille and visiting family in Oslo over Christmas and dropped in at the University where we briefly met. Later I got to know him very well, not only through our work and joint papers (15 altogether) but through holidays taken together particularly skiing ones in Rondane and Sjusjøen. He spent 6 months with me in Warwick in 1982, funded by an EPSRC visiting fellowship, which was the start of our collaboration, followed by many months together at further research visits not only at Warwick, e.g. the symposium in 1986-87 on *Operator Algebras and Applications* but also in Trondheim and Oslo, Swansea and Cardiff, Kyoto and Sapporo, Canberra and Bangkok, Ottawa and Toronto, Mittag-Leffler and Fields, Waterloo etc spread over many years.

Our first conference together was at Arco Felice, Naples in March 1978. Ola and Aki drove there from Marseille in his Citroën deux chevaux. This was a meeting on *Mathematical Problems in Quantum Theory of Irreversible Processes* organized by Vittorio Gorini - which brought together our mutual interests in derivations and generators of one parameter semigroups of positive maps - on which we would later collaborate. On the day after the conference Ola drove Aki, John Lewis, Geoffrey Sewell and myself to Pompeii and the Bay of Naples - all in his deux chevaux! At a meeting in Chennai organized by Sunder, Ola took me on an expedition to find saffron in the local markets, which was later put to good use in his bouillabaisse, the finest I have ever had, with fresh fish from the harbor back in Norway.

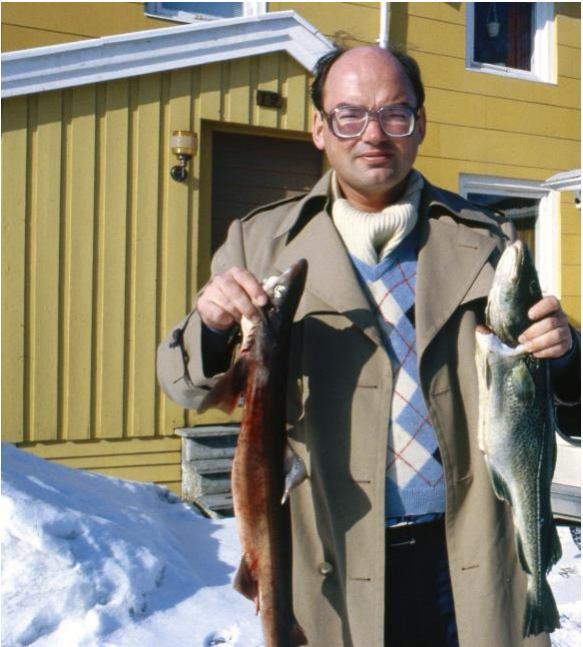
Ola preparing Bouillabaisse in Trondheim.

Ola classified in his Phd thesis [1], under the supervision of Erling Størmer at Oslo, AF C*-algebras, inductive limits of finite dimensional operator algebras in terms of what are now known as Bratteli diagrams. This work was not only pivotal in Elliott's classification of AF-algebras through K-theory; *J of Algebra* 1976, but are ubiquitous in operator algebras and dynamical systems. In the theory of subfactors, the diagrams from the towers of relative commutants encode the principal graph and dual principal graphs (Jones, Ocneanu, Popa *et al*). In Elliott *J of Algebra* 1976, the Bratteli diagram encodes the K-theoretic data of the morphisms between approximating finite dimensional algebras and understanding intertwining pairs of Bratelli diagrams leads to the classification. In subfactor theory with inclusions of a pair of towers of relative commutants, there is a further degree of freedom measured as a connection in the language of Ocneanu. The space of infinite paths in Bratteli diagrams and associated odometers and dynamical systems are fundamental to the classification of symmetries on Cantor sets. Indeed every minimal Z-action on a Cantor set is equivalent to a minimal Bratteli-Vershik system; Herman, Putnam and Skau, *Inter J of Math* 1992.

My own work with Ola started with dynamical systems of one parameter semigroups of positive maps, begun during his long term visit to Warwick in 1982; Bratteli and Evans, *Ergodic Thy & Dyn. Sys.* 1983. Later visits to Swansea in the summers of 1998, 1999 and 1990, led to the collaboration with George Elliott and Akitaka Kishimoto on non commutative spheres, the irrational rotation algebra and the K-theoretic obstructions to classifying such amenable C*-algebras and their dynamical systems. This was particularly motivated by the remarkable construction by Blackadar, *Annals of Math* 1990 of a $\mathbb{Z}_2$-symmetry on the Fermion algebra with non AF fixed point algebra by building the Fermion algebra in a non-standard way as an inductive limit of matrix algebras over the circle: $C(\mathbb{T}, M_{4^n})$ and the subsequent work of Kumjian in showing that a $\mathbb{Z}_2$-symmetry of a Bunce-Deddens algebra and the corresponding crossed product $C(\mathbb{T}) \rtimes$ dyadic rotations $\rtimes \mathbb{Z}_2$ yielded an AF algebra. The issue of existence of such symmetries had been a well known open problem – c.f. the situation with von Neumann algebras where there is an unique outer action of $\mathbb{Z}_2$ on the hyperfinite $II_1$ factor factor and indeed there is an unique hyperfinite subfactor of index 2.

This led us to consider the $\mathbb{Z}_2$ symmetry on the rotation algebra with the *flip* σ; $\sigma(u) = u^{-1}, \sigma(v) = v^{-1}$ on the generators. On the classical two torus, $\mathbb{T}^2$ this yields a singular orbifold $\mathbb{T}^2/\mathbb{Z}_2$, which can be thought of as a tetrahedron, topologically a sphere, but with four singular vertices. This is better studied through the crossed product $C(\mathbb{T}^2) \rtimes \mathbb{Z}_2$ as $M_2$ valued functions on the sphere or tetrahedron, restricted to

$\mathbb{C}^2$ at four singular points, the vertices of the tetrahedron. This led us to call the fixed point algebras and crossed products of rotation algebras as non commutative spheres but as Alain Connes pointed out, they are better described as non commutative toroidal orbifolds as they do not have the *K*-theory of a sphere. The origin of the naming of the term *orbifold* is described by Thurston on page 300 of Chapter 13 of his book, The *Geometry and Topology of three-manifolds,* 1978. Blackadar's construction of a symmetry of the Fermion algebra *A* with non-AF fixed point algebra then drew us to study orbifolds of the circle, or an interval:

Taking matrix valued functions constrained at the singular points to have dimension drops and their inductive limits led us on a path towards:

- *G* finite on UHF $A = \otimes_N M_{|G|}$ with $A^G$ not AF; Bratteli, Elliott, Evans and Kishimoto, *K-theory* 1994,
- *G* compact on UHF with non AF fixed point algebra; Evans and Kishimoto, *J Funct. Anal* 1991,
- $C(\text{Cantor}) \rtimes (\mathbb{Z} \rtimes \mathbb{Z}_2)$ is AF when there is a fixed point for the $\mathbb{Z}_2$ action and the $\mathbb{Z}$ action is minimal; Bratteli, Evans and Kishimoto, *Ergodic Thy. and Dyn. Sys.*1993,
- there exist non AF C*-algebras which when tensored with a certain UHF algebra become UHF; Evans and Kishimoto, *J Funct. Anal* 1991,
- the fixed point algebra of an irrational rotation algebra by the flip is AF; Bratteli and Kishimoto, *Comm Math Phys*. 1992,
- the irrational rotation algebras are inductive limits of sums of matrix algebras over the continuous functions on a circle; Elliott and Evans, *Annals of Math*, 1993.

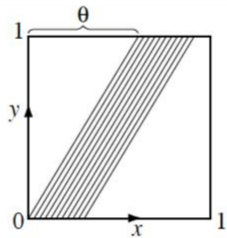

The flip on the irrational or Kronecker flow on the torus yields a non-singular flow on the sphere which is zero dimensional in a strong sense – the corresponding C*-algebra is approximately finite dimensional.

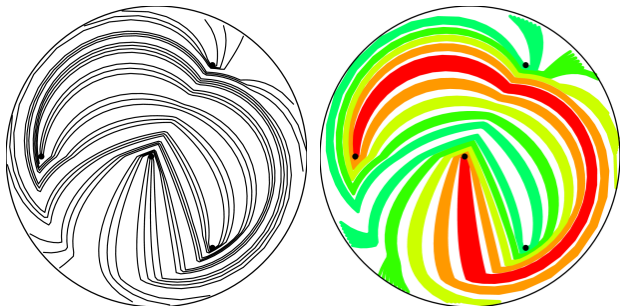

Our program was also influenced by the work of Putnam, *Pacific J Math* 1989, who showed that the crossed product $C(\text{Cantor}) \rtimes \mathbb{Z}$ of a minimal action on a Cantor set is AT, i.e. an inductive limit of sums of matrix algebras over the continuous functions on a circle, and Dadarlat and Nemethi, *J of Operator Thy.* 1990, who classified finite dimensional algebras over CW complexes by $K_0$ and $K_1$ up to shape equivalence. Here, the approximating intertwining Bratteli diagrams of morphisms are valid up to homotopy.

This work described above laid the foundations for subsequent work over the last 25 years on the classification of amenable C*-algebras by *K*-theoretic data - before which the classification was out of sight and did not appear feasible.

Our work on orbifolds also directly lead to the study of orbifold subfactors, Evans and Kawahigashi, *Comm Math Phys* 1994, as reported in the first Danish-Norwegian Workshop on Operator Algebras in Røros in 1991. Just as had been understood in the C*-setting initiated by Blackadar, symmetries on subfactors, through symmetries on principal graphs, dual principal graphs and connections with no fixed points could yield new subfactors as orbifolds or simultaneous fixed point algebras $N^G \subset M^G$ with the same index as the original subfactor $N \subset M$ but different principal and principal graphs. For example, for index less than 4, the oribifolds of the $A_{2n-3}$ subfactors yield the $D_n$ subfactors. In the case *n*=3 of the Dynkin diagram $A_3$ which is self dual with orbifold $D_3 = A_3$, this is related to high temperature low temperature duality of Krammers-Wannier in the Ising model - an orbifold. Generalizations in the Potts model with underlying symmetry or orbifolds appear in the understanding of the double of the Haagerup subfactor as the grafting of two orbifolds of two Potts models in current work of Evans and Gannon.

Ola had a long connection and affection for Thailand making regular visits – his wife for more than 30 years Wasana was from Phitsanulok. In 2016, Paulo Bertozzini initiated at Thammasat University in Bangkok, the *Ola Bratteli Mathematical Physics and Mathematics in Thailand Colloquium* when I was honoured to give the first talk.

Ola had a generous spirit and integrity. We will miss his presence and friendship.

### Acknowledgements

Editing has been done by Magnus B. Landstad, translations from Norwegian by Margaret Forbes, and Bratteli diagrams by Petter Nyland.